\documentclass[smallextended]{svjour3}       
\usepackage{latexsym}
\usepackage{amssymb}
\usepackage{amsmath}
\usepackage[mathscr]{eucal}
\usepackage{graphicx}
\usepackage{hyperref}
\usepackage{caption}
\usepackage{subcaption}

\renewcommand{\qed}{\hfill{\ \ \rule{2mm}{2mm}} \vspace{0.2in}}

\newcommand{\ind}{1\hspace{-2.3mm}{1}}

\begin{document}



\title{Deviation Estimates for Extremal Relay Random Geometric Graphs}
\author{ \textbf{Ghurumuruhan Ganesan}
\thanks{E-Mail: \texttt{gganesan82@gmail.com} } \\
\ \\
IISER, Bhopal}
\date{}
\maketitle

\begin{abstract}
In this paper, we consider a deterministic graph~\(\Gamma\) drawn on the unit square with straight line segments as edges and connect vertices of~\(\Gamma\) using edges of a random geometric graph (RGG)~\(G\) with adjacency distance~\(r_n\) as relays. We call the resulting graph as a \emph{relay} RGG and determine sufficient conditions under such relay RGGs exist and are also near optimal, in terms of the graph parameters of~\(\Gamma.\) We then equip edges of~\(G\) with independent, exponentially distributed weights and obtain bounds for the maximum possible weight~\(W_n\) of a relay RGG with a given length~\(L_n.\)

\vspace{0.1in} \noindent \textbf{Key words:} Relay Random Geometric Graphs; Minimum Length; Maximum Weight; Deviation Estimates.

\vspace{0.1in} \noindent \textbf{AMS 2000 Subject Classification:} Primary: 60C05, 60D05;
\end{abstract}

\bigskip

\renewcommand{\theequation}{\thesection.\arabic{equation}}
\setcounter{equation}{0}
\section{Introduction} \label{intro}

Random geometric graphs (RGGs) are a class of Euclidean random graphs that often occur in the study of wireless networks~\cite{gupta}~\cite{diaz}~\cite{penrose}. Consequently, structural properties of RGGs have been studied extensively from both theoretical and application perspectives. For e.g., in~\cite{gan_one}~\cite{gan_two} we have obtained bounds on the size of the giant component and diameter in subcritical RGGs, respectively, using percolation theoretic techniques. As examples of applications,~\cite{zhang}\cite{song}\cite{ferrero} use RGGs for analysis of different models of vehicular and wireless networks. Many variants of RGGs have been studied as well: In~\cite{david}, random ``anti"-geometric graphs are defined and connectivity properties are studied for the large vertex regime and recently,~\cite{pen2} studied ``soft" RGGs where edges of RGG are removed with a certain probability and connectivity is again established for various regimes of the graph parameters.

In this paper, we consider a weighted subclass of RGGs which we call as relay RGGs and obtain deviation estimates for its extremal counterparts. Specifically, we begin with a deterministic planar graph~\(\Gamma\) with edges drawn as straight line segments in the unit square and use the edges of a random geometric graph (RGG) as ``relays" to connect vertices of~\(\Gamma\) using \emph{edge disjoint} paths (if such paths exist). We call the resulting graph as a relay RGG and establish conditions under which~\(\Gamma\) is ``supported" in the sense that relay RGGs exist and are near optimal in terms of length.

We then equip each edge  of~\(G\) with  a random weight whose entries are independent and identically distributed (i.i.d.) and obtain deviations for the  \emph{maximum weight} of relay RGGs.  To motivate this, consider a typical wireless application scenario where the vertices of~\(\Gamma\) represent communication terminals and the edges of the RGGs are relay links through which messages are passed from one vertex of~\(\Gamma\) to another, without interference (and hence the edge disjoint condition).

When the edges of the relay RGG act as communication links, the  weights have a natural interpretation of the \emph{channel fading} gain and it is of interest to maximize the network throughput by choosing the relay RGG with the maximum possible weight. Selecting links with maximum fading gain (also known as diversity gain) has important applications in wireless communications and for more details on exploiting diversity in wireless networks, we refer to Chapter~\(7,\)~\cite{goldsmith}, Chapter~\(12,\)~\cite{molisch} and references therein. In this paper we obtain deviation estimates for the maximum weight of relay RGGs that are of the same order and also establish sufficient conditions for the~\(L^2-\)convergence of the maximum weight, appropriately scaled and centred.



In the following section, we state and prove our two main results regarding minimum length and maximum weight relay RGGs.

\setcounter{equation}{0}
\renewcommand\theequation{\thesection.\arabic{equation}}
\section{Relay Random Geometric Graphs}\label{sec_ramsey}
Let~\(f\) be any density on the unit square~\(S = \left[-\frac{1}{2},\frac{1}{2}\right]^2\) satisfying
\begin{equation}\label{f_eq}
\epsilon_1 \leq \inf_{x \in S} f(x) \leq \sup_{x \in S} f(x) \leq \epsilon_2
\end{equation}
for some constants~\(0 < \epsilon_1 \leq \epsilon_2 < \infty.\) We assume throughout that constants do not depend on~\(n.\)

For~\(n \geq 1,\) let~\(\{X_i\}_{1 \leq  i \leq n}\) be independently and identically distributed (i.i.d.) random variables with density~\(f(.).\) We assign the label~\(i\) to the vertex location~\(X_i\) and let~\(K_n\) be the complete graph with vertex set~\(\{1,2,\ldots, n\}.\)  For~\(r_n > 0,\) let~\(G  = G(r_n)\subset K_n\) be the subgraph of~\(K_n\) consisting of all edges of length at most~\(r_n.\) Here and henceforth length always refers to the Euclidean length and we define~\(G\) to be the random geometric graph (RGG) with adjacency distance~\(r_n.\)

Let~\(\Gamma = (V(\Gamma),E(\Gamma))\) be any deterministic graph containing~\(v_0 = v_0(n)\) vertices and~\(e_0 = e_0(n)\) edges in the interior of the unit square~\(S.\) Let~\(G_{loc} = G_{loc}(G,\Gamma)\) be the graph whose edge set is the union of the edge set of~\(G\) and the edge of the graph obtained by connecting each vertex~\(v \in V(\Gamma)\) to all vertices of~\(G\) within distance~\(r_n\) from~\(v.\)

A path of length~\(l\) in~\(G_{loc}\) is a sequence~\(\pi = (a_1,\ldots,a_{l+1})\) containing~\(l\) edges, where each vertex~\(a_i\) is adjacent to~\(a_{i+1}.\) The vertices~\(a_1\) and~\(a_{l+1}\) are defined to be the endvertices of~\(\pi.\) We say~\(\pi\) is a \emph{relay} path if~\(a_1\) and~\(a_{l+1}\) are vertices of~\(\Gamma\) and every other~\(a_i\) is a vertex of~\(G.\) Correspondingly, we say that a subgraph~\(H \subset G_{loc}\) is a \emph{relay} RGG (with respect to~\(\Gamma\)) if~\(H= \bigcup_{f \in E(\Gamma)} P(f)\) where:\\
\((i)\) Each~\(P(f)\) is a relay path with the same endvertices as~\(f,\)\\
\((ii)\) The path~\(P(f_1)\) is vertex disjoint from~\(P(f_2)\) if~\(f_1 \) and~\(f_2\) do not share an endvertex and\\
\((iii)\) If~\(f_1\) and~\(f_2\) share a common endvertex~\(u,\) then~\(P(f_1)\) and~\(P(f_2)\) have only the vertex~\(u\) in common.

Intuitively, if the number of edges in~\(\Gamma\) is too large, then there may not exist a relay RGG with respect to~\(\Gamma.\) To illustrate, consider the graph~\(\Gamma_0\) formed by the union of the edges~\(f_i = (u_i,v_i), 1 \leq i \leq n\) where~\(u_i = \left(\frac{i-1}{n},0\right)\) and~\(v_i = \left(\frac{i-1}{n},\frac{2}{n}\right).\) Thus the edges~\(\{f_i\}\) are vertical, parallel and spaced~\(\frac{1}{n}\) apart and each~\(f_i\) has length~\(\frac{2}{n}.\) If the adjacency distance~\(r_n = \frac{1}{2n},\) then we see that there does not exist a relay RGG with respect to~\(\Gamma.\)


Let~\(L(G,\Gamma)\) be the minimum number of edges in a relay RGG, if it exists; else set~\(L(G,\Gamma) = \infty.\) The following result determines sufficient conditions under which relay RGGs exist and are also near optimal, in terms of the overall number of edges. For an edge~\(f \in E(\Gamma)\) we let~\(l(f)\) be the Euclidean length of~\(f\) and define~\(l_0 := \min_{f} l(f)\) and~\(l_{tot} := \sum_{f} l(f)\) respectively to be the minimum length of an edge and the sum total of all edge lengths, in~\(\Gamma.\) Throughout constants do not depend on~\(n.\)
\begin{theorem}\label{thm_edge_har}  Suppose~\(\Gamma\) has~\(e_0\) edges and
\begin{equation}\label{gam_cond}
e_0 \leq n^{\alpha},\;\;\;\frac{1}{n^{\beta}} \leq r_n \longrightarrow 0\;\;\text{ and }\;\;\frac{l_0}{r_n} \longrightarrow \infty
\end{equation}
where~\(0 < \alpha,\beta < 1\) are constants satisfying~\(\beta < \frac{1-\alpha}{2}\) so that~\(e_0\) is much smaller than~\(nr_n^2.\) for every~\(\epsilon > 0\) there is a constant~\(C  >0\) such that
\begin{equation}\label{har_weak}
\mathbb{P}\left( \frac{l_{tot}}{r_n} \leq L(G,\Gamma) \leq \frac{l_{tot}}{r_n}(1+\epsilon)\right) \geq 1- e^{-C nr_n^2}.
\end{equation}
Moreover, if~\(0 < \beta < \frac{1-\alpha}{4},\) then there is a constant~\(D > 0\) such that
\begin{equation}\label{har_str}
\mathbb{P}\left( \frac{l_{tot}}{r_n}  \leq L(G,\Gamma) \leq \frac{l_{tot}}{r_n} +2e_0\right) \geq 1- e^{-D nr_n^4}.
\end{equation}
\end{theorem}
As a check we see that~\(l_{tot} \geq l_0 e_0\) and so the term~\(\frac{l_{tot}}{r_n} \geq \frac{l_0}{r_n} \cdot e_0\) is much larger than~\(e_0,\) by~(\ref{gam_cond}). In other words, as long as the number of vertices in~\(\Gamma\) is small enough in the sense of~(\ref{gam_cond}), then we can find a relay RGG using the edges of~\(G\) that is also near optimal in terms of number of edges.





To prove Theorem~\ref{thm_edge_har}, we use the following result of independent interest that obtains  two point concentration for the number of edges in the shortest length relay path of an edge~\(f = (u,v) \in \Gamma\) in the graph~\(G_{loc}.\) If~\(d(u,v)\) and~\(d_{gr}(u,v)\) respectively denote the Euclidean and graph distance in~\(G_{loc}\) between vertices~\(u\) and~\(v,\) then~\(d_{gr}(u,v) > \frac{d(u,v)}{r_n}\) since each edge in~\(G_{loc}\) has Euclidean length strictly less than~\(r_n.\) For~\(\epsilon > 0,\) we define the events~\[E_{u,v} = E_{u,v}(\epsilon):= \left\{ d_{gr}(u,v) \leq \frac{d(u,v)}{r_n}(1+\epsilon)\right\}\]  and \[F_{u,v} := \{d_{gr}(u,v) \in \{d_{u,v},d_{u,v}+1\}\},\] where~\(d_{u,v}\) is the smallest integer strictly larger than~\(\frac{d(u,v)}{r_n}\) and have the following result.
\begin{lemma}\label{thm_dist_rgg}  For every~\(\epsilon > 0\) there are constants~\(M,D > 0\) not depending on the choice of~\(u\) or~\(v\) such that if~\(r_n \geq \sqrt{\frac{M \log{n}}{n}}\) then
\begin{equation}\label{eps_bound}
\mathbb{P}(E_{u,v}) \geq 1 - e^{-Dnr_n^2}
\end{equation}
and
\begin{equation}\label{two_point}
\mathbb{P}\left(F_{u,v}\right) \geq 1- \frac{2d(u,v)}{r_n} \cdot \exp\left(-\frac{D nr_n^4}{d^2(u,v)}\right).
\end{equation}
\end{lemma}
For example, if~\(d(u,v) \leq \frac{r_n^2 \sqrt{n}}{\log{n}}\) then we get from~(\ref{two_point}) that~\(F_{u,v}\) occurs with probability
\begin{equation}\label{two_point_eval}
\mathbb{P}(F_{u,v}) \geq 1-\frac{2d(u,v)}{r_n} \cdot e^{-D(\log{n})^2} \geq 1- 2\sqrt{2n}\cdot e^{-D(\log{n})^2}
\end{equation}
since~\(d(u,v) \leq \sqrt{2}\) and~\(\frac{1}{r_n} \leq \sqrt{n}\) by Lemma statement.


\emph{Proof of Lemma~\ref{thm_dist_rgg}}:  For simplicity we set~\(r_n = r\) and assume that~\(u=(0,0)\) and~\(v=(Kr,0)\) where~\(K := \frac{d(u,v)}{r} > 1\) and compute their graph distance in~\(G_{loc}.\) For real numbers~\(0<\eta_1 < \eta_2\) to be determined later, let~\(\delta \in (\eta_1,\eta_2)\) be such that~\(W := \frac{K}{1-\delta}\) is an integer. This is possible if~\(\eta_1\) and~\(\eta_2\) are such that
\begin{equation}\label{eta_cond_one}
\frac{K}{(1-\eta_2)} - \frac{K}{(1-\eta_1)} = \frac{K(\eta_2-\eta_1)}{(1-\eta_1)(1-\eta_2)} \geq 1.
\end{equation}

\begin{figure}[tbp]
\centering
\includegraphics[width=4in, trim= 40 550 70 90, clip=true]{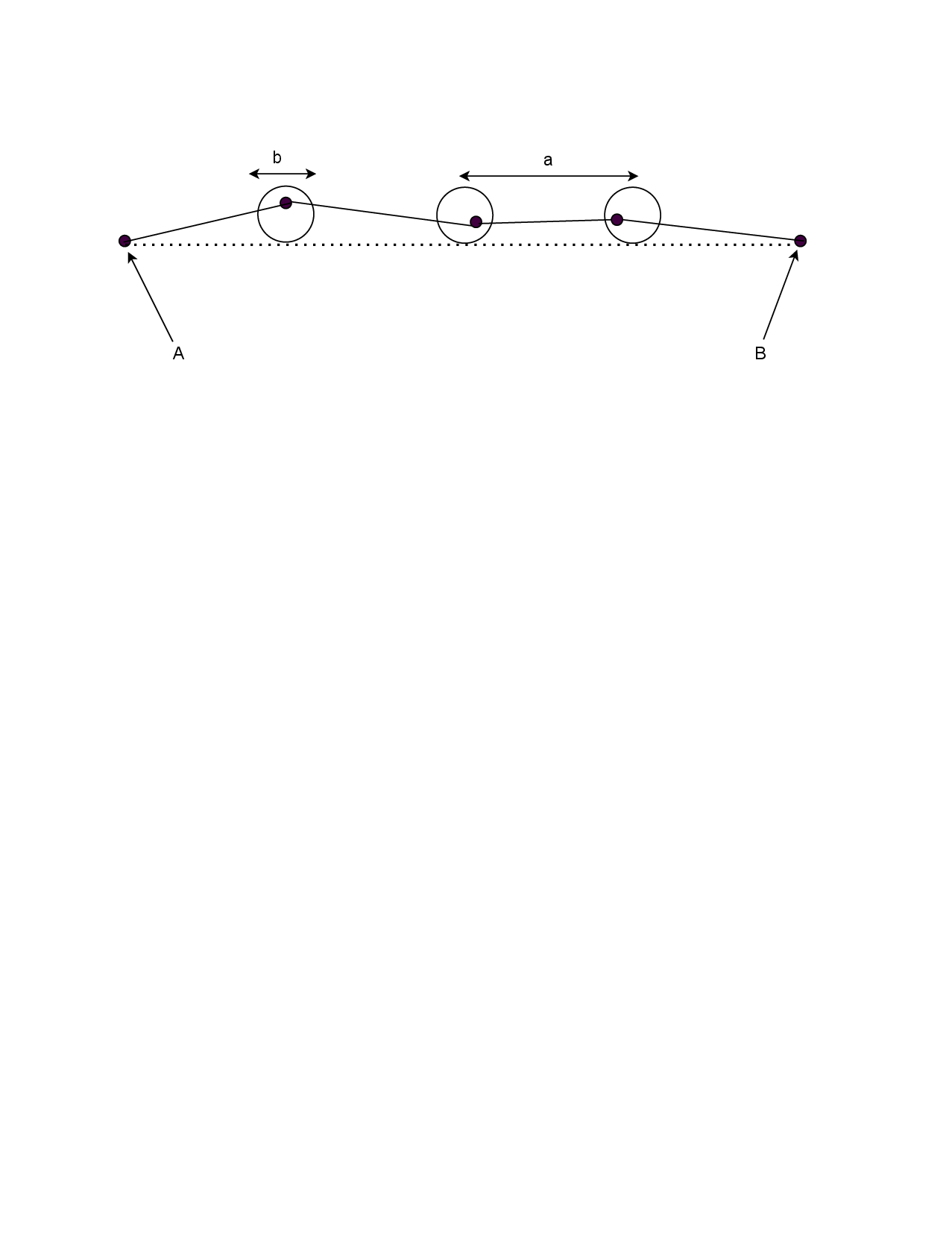}
\caption{The circles~\(\{S_i\}_{1 \leq i \leq W-1},\) each of diameter~\(b = \frac{2\delta r}{L} ,\) with corresponding centres that are spaced~\(a = r(1-\delta)\) apart. The vertices~\(u\) and~\(v\) are represented by~\(A = (0,0)\) and~\(B = (Kr,0),\) respectively and the dark dots represent the vertices of~\(G_{loc}\) that form the path~\(\pi\) from~\(A\) to~\(B,\) denoted by the solid line segments. }
\label{fig_1}
\end{figure}

For~\(1 \leq i \leq W=1\) let~\(S_i\) be the circle with centre~\(\left(i\cdot r\cdot(1-\delta),\frac{\delta r}{L}\right)\) and radius~\(\frac{\delta r}{L}\) (see Figure~\ref{fig_1}) where~\(L > 1\) is a constant to be determined later. For generality, we have chosen these circle centres to also account for the case where~\(u\) or~\(v\) could be close to or on the boundary of the unit square. For~\(1 \leq l \leq W-1,\) we have that~\[\mathbb{P}(X_j \in S_i) = \int_{S_i} f(z)dz \geq \frac{\epsilon_1 \pi \delta^2 r^2}{L^2},\] by~(\ref{f_eq}). Defining~\(E_i\) to be the event that~\(S_i\) is nonempty we have
\begin{equation}
\mathbb{P}(E_i) \geq 1-\left(1-\frac{\pi \epsilon_1 \delta^2 r^2}{L^2}\right)^{n}  \geq 1-\exp\left(-\frac{\pi\epsilon_1 \delta^2 nr^2}{L^2}\right)  \nonumber
\end{equation}
using~\(1-y \leq e^{-y}\) for~\(y > 0.\) If~\(E_{tot} = \bigcap_{i=1}^{W}E_i,\) then an application of the union bound gives
\begin{equation}
\mathbb{P}\left(E_{tot}\right) \geq 1 - W\cdot \exp\left(-\frac{\pi \epsilon_1 \delta^2 nr^2}{L^2}\right) \geq 1- \frac{K}{1-\eta_2}\exp\left(-\frac{\pi \epsilon_1 \eta_1^2 nr^2}{L^2}\right),\\
\label{e_occ_est}
\end{equation}
since~\(\delta > \eta_1\) and~\(W = \frac{K}{1-\delta} \leq \frac{K}{1-\eta_2}.\)


Suppose~\(E_{tot}\) occurs and for~\(1 \leq  l \leq W-1\) pick a vertex~\(v_l \in S_l\) and set~\(v_0 = (0,0)\) and~\(v_W = (Kr,0).\) By construction, the Euclidean distance between~\(v_i\) and~\(v_{i+1}\) is at most~\(r(1-\delta) + \frac{2\delta r}{L} < r\) provided~\(L > 2.\) Similarly the distance between~\(v_i\) and~\(v_{i+1}\) is at least~\(r\left(1-\delta - \frac{2\delta}{L}\right) \geq r\left(1-\frac{3\delta}{2}\right),\) provided~\(L > 4.\)
Finally, we see that each circle in Figure~\ref{fig_1} has diameter at most~\(\frac{2\delta r}{L} \leq \frac{\delta r}{8} < \frac{r}{8},\) again provided~\(L > 16.\) Fixing such an~\(L,\) we then see that the path~\(\pi = (v_0,v_1,\ldots,v_W)\) satisfies the following properties:\\
\((p1)\) There are~\(W = \frac{K}{1-\delta}\) edges in~\(\pi,\)  \\
\((p2)\) Any two vertices in~\(\pi\) are at a distance of at least~\(r\left(1-\frac{3\delta}{2}\right)\) and\\
\((p3)\) The graph distance between~\(v_0\) and~\(v_W\) satisfies
\begin{equation}\label{dg_est}
d_{gr}(v_0,v_W) - K \leq \frac{K}{1-\delta} - K =  \frac{K\delta}{1-\delta} \leq 1+\gamma
\end{equation}
provided~\(\delta < \eta_2 \leq \frac{1+\gamma}{1+\gamma+K}.\) We fix such a~\(\delta\) henceforth.

Summarizing, we need to choose~\(L>16\) and~\(0 < \eta_1 < \eta_2 \leq \frac{1+\gamma}{1+\gamma+K}\)
in such  a way that~(\ref{eta_cond_one}) also holds. Setting~\(\eta_2 = \frac{1+\gamma}{1+\gamma+K},\) the condition~(\ref{eta_cond_one}) reduces to~\(\eta_1 \leq \frac{\gamma}{K+\gamma}.\) Therefore we fix~\(\eta_1 = \frac{\gamma}{K+\gamma}\) and get from~(\ref{e_occ_est}) that
\begin{equation}\label{e_occ_est2}
\mathbb{P}(E_{occ}) \geq 1-(1+K+\gamma)\cdot \exp\left(-\frac{\pi\epsilon_1 nr^2}{L^2}\left(\frac{\gamma}{K+\gamma}\right)^2\right).
\end{equation}

Suppose the event~\(E_{occ}\) occurs. We then know from~(\ref{dg_est}) that~\(d_G(v_0,v_W) \leq K+\gamma+1\) and so setting~\(\gamma = K\epsilon\) in~(\ref{e_occ_est2}) and defining~\(D_1 = \frac{\pi\epsilon_1}{L^2},\) then we get that
\begin{eqnarray}
\mathbb{P}\left(d_G(v_0,v_W) \leq K(1+\epsilon)+1\right) &\geq& 1- (1+K(1+\epsilon))\exp\left(-D_1nr_n^2\frac{\epsilon^2}{(1+\epsilon)^2}\right) \nonumber\\
&\geq& 1- \frac{2}{r_n} \exp\left(-D_2 nr_n^2\right) \nonumber\\
&\geq& 1- \exp\left(-D_3nr_n^2\right) \nonumber
\end{eqnarray}
for some constants~\(D_2,D_3 > 0,\) provided~\(r_n \geq \sqrt{\frac{M\log{n}}{n}}\) and~\(M\) is a large enough constant. Fixing such an~\(M,\) we get~(\ref{eps_bound}).

To obtain~(\ref{two_point}), we choose~\(\gamma = \frac{1}{2}\) in~(\ref{e_occ_est2}), so that~\[K+ \gamma+1  = K + \frac{3}{2} < \alpha(K)+2,\] strictly, where~\(\alpha(K)\) is the smallest integer strictly larger than~\(K.\) Recalling the definition of the event~\(F_{u,v}\) prior to the statement of Lemma~\ref{thm_dist_rgg}, we get from~(\ref{e_occ_est2}) that
\begin{equation}\label{damini_tits}
\mathbb{P}(F_{u,v}) \geq 1-2K\exp\left(-C_1\frac{nr_n^2}{K^2}\right).
\end{equation}
Substituting~\(K = \frac{d(u,v)}{r},\) we get~(\ref{two_point}) and this completes the proof of the Lemma.~\(\qed\)

The final ingredient in our proof of Theorem~\ref{thm_edge_har} is the following deviation estimate regarding sums of independent  Bernoulli random variables.  Let~\(Y_1,\ldots,Y_m\) be independent Bernoulli random variables satisfying~\[\mathbb{P}(Y_j = 1) =1-\mathbb{P}(Y_j =0) >0.\] If~\(T_m := \sum_{j=1}^{m} Y_j\) and~\(\mu_m = \mathbb{E}T_m,\) then for any~\(0 < \epsilon \leq \frac{1}{2}\) we have that
\begin{equation}\label{conc_est_f}
\mathbb{P}\left(|T_m - \mu_m| \geq \epsilon \mu_m\right) \leq \exp\left(-\frac{\epsilon^2}{4} \mu_m\right).
\end{equation}
For a proof of~(\ref{conc_est_f}), we refer to Corollary A.1.14, pp. 312 of Alon and Spencer (2008).

\emph{Proof of Theorem~\ref{thm_edge_har}}: Let~\(f =(u,v) \in \Gamma\) be an edge of length~\(l(f).\) In any relay RGG~\(G_{\Gamma}\) the path~\(P(f)\) containing~\(u\) and~\(v\) as endvertices has at least~\(\frac{l(f)}{r_n}\) edges. Summing over all the edges~\(f\) of~\(\Gamma\) gives the lower bound in~(\ref{har_weak}) and in~(\ref{har_str}).

For the upper bound in~(\ref{har_str}), we argue as follows. Suppose initially that~\(\Gamma\) consists of a single edge~\(f = (u,v)\) with~\(u=(0,0)\) and~\(v = (Kr_n,0)\) where~\(K := \frac{l(f)}{r_n}\) satisfies
\begin{equation}\label{k_bounds}
\frac{2}{r_n} \geq K \geq \frac{l_0}{r_n} \longrightarrow \infty,
\end{equation}
by~(\ref{gam_cond}). We consider equally spaced circles~\(\{S_i(f)\}_{1 \leq i \leq N(f)}\) as in Figure~\ref{fig_1} with~\(A = u\) and~\(B=v\) and let~\(R_i(f)\) be the number of vertices of~\(\{X_i\}\) present in~\(S_i(f).\) Recalling the definition of~\(\delta\) in the proof of Lemma~\ref{thm_dist_rgg}, we have that~\(\mathbb{E}R_i(f) \geq  2C_1 \delta^2nr_n^2\) for some constant~\(C_1 >0.\) All constants henceforth do not depend on the choice of~\(f.\) Defining the event~\[E_i(f) := \{R_i(f) \geq C_1 n\delta^2 r_n^2\},\] we therefore get from~(\ref{conc_est_f}) that
\begin{equation}\label{damini_tits_ax}
\mathbb{P}(E_i(f)) \geq 1-e^{-C_2n\delta^2r_n^2}
\end{equation}
for some constant~\(C_2 >0.\)

As in the proof of~(\ref{damini_tits}), we choose~\(\gamma = \frac{1}{2}\) so that
\begin{equation}\label{del_rel}
\delta  > \eta_1 = \frac{\gamma}{\gamma+K} \geq  \frac{\gamma}{2K} \geq C_3r_n
\end{equation}
for some constant~\(C_3 > 0,\) where the second inequality in~(\ref{del_rel}) follows from the fact that~\(K \geq \frac{l_0}{r_n} \longrightarrow \infty\) (see Theorem statement) and the final relation in~(\ref{del_rel}) is true since~\(K \leq \frac{2}{r_n}\) by definition. For future use we also remark that as in the proof of Lemma~\ref{thm_dist_rgg} the parameter~\(\gamma \in \left\{\frac{1}{2},K\epsilon\right\}\) with~\(\epsilon  >0\) as in Theorem statement and so
\begin{equation}\label{del_rel_two}
\delta \leq \eta_2 = \frac{1+\gamma}{1+\gamma+K} \leq \frac{1+K\epsilon}{1+K(1+\epsilon)} \leq 2\epsilon \leq \frac{1}{4}
\end{equation}
for all~\(n\) large provided~\(\epsilon > 0\) is small. Plugging the final estimate of~(\ref{del_rel}) into~(\ref{damini_tits_ax}), we get
\begin{equation}\label{damini_tits_ax_2}
\mathbb{P}(E_i(f)) \geq 1-e^{-4C_4nr_n^4}
\end{equation}
for some constant~\(C_4 >0.\)

From the discussion prior to~(\ref{e_occ_est2}), we see that the number of circles in consideration is~\[N(f) \leq \frac{K}{1-\eta_2} \leq 1+K+\gamma \leq 2K\] by the lower bound in~(\ref{k_bounds}). Thus defining~\(E_{tot}(f) := \bigcap_{i=1}^{N(f)} E_i(f),\) we get from~(\ref{damini_tits_ax_2}), the upper bound in~(\ref{k_bounds}) and the union bound that
\begin{equation}
\mathbb{P}(E_{tot}(f)) \geq 1- 2K\cdot e^{-4C_4nr_n^4}\geq 1- \frac{4}{r_n} \cdot e^{-4C_4nr_n^4} \geq 1-e^{-2C_4nr_n^4} \label{e_tot_est_ax3}
\end{equation}
since~\(r_n \geq \frac{1}{n^{\beta}}\) with~\(\beta < \frac{1}{4}\) strictly.

Suppose~\(E_{tot}(f)\) occurs. Arguing as in the proof of Lemma~\ref{thm_dist_rgg}, we see that there exists a relay path~\(P(f)\) containing~\(u\) and~\(v\) as endvertices, whose length is either~\(\alpha_n(f)\) or~\(\alpha_n(f)~+~1,\) where~\(\alpha_n(f)\) is the smallest integer strictly larger than~\(\frac{l(f)}{r_n}.\) This obtains the upper bound in~(\ref{har_str}) for the case when~\(\Gamma\) consists of a single edge.

Next, if~\(\Gamma\) consists of multiple edges, then we label the edges of~\(\Gamma\) as~\(h_1,\ldots,h_{e_0}\) and perform the same analysis as above, iteratively. Letting~\(E_{edge} := \bigcap_{i=1}^{e_0} E_{tot}(h_i),\) we get from~(\ref{e_tot_est_ax3}) and the union bound that
\begin{equation}\label{e_edge_est}
\mathbb{P}(E_{edge}) \geq 1-e_0\cdot e^{-2C_5nr_n^4} \geq 1-e^{-C_5nr_n^4},
\end{equation} since~\(e_0 \leq n^{\alpha}.\) Assuming~\(E_{edge}\) occurs, our strategy now is to pick the desired relay paths~\(\{P(h_i)\}\) iteratively in such a way that~\(P(h_i)\) and~\(P(h_j), j < i\) share no vertex in common if~\(h_i\) and~\(h_j\) do not share an endvertex and share the common endvertex of~\(h_i\) and~\(h_j,\) otherwise.

Indeed, consider the set of circles~\(S_1(h_i),\ldots,S_W(h_i)\) as in Figure~\ref{fig_1} used to construct the relay path~\(P(h_i).\) From  the occurrence of~\(E_{edge},\) we see that each circle~\(S_k(h_i)\) contains at least~\(Cnr_n^4 \geq Cn^{1-4\beta}\) vertices of~\(\{X_l\}.\) Moreover, based on the discussion prior to properties~\((p1)-(p3)\) in the proof of Lemma~\ref{thm_dist_rgg}, we know that~\(S_k(h_i)\) has diameter at most~\(\frac{r_n}{8}\) and from property~\((p2)\) we also know that any two vertices in~\(P(h_j),j <i\) are spaced at least~\(r_n\left(1-\frac{3\delta}{2}\right) \geq \frac{5r_n}{8}\) apart, by~(\ref{del_rel_two}). Thus at most one vertex of~\(P(h_j)\) is present in~\(S_k(h_i)\) and so at most
\begin{equation}\label{e_not_est}
e_0 \leq n^{\alpha} = o(n^{1-4\beta})
\end{equation}
vertices of~\(S_k(h_i)\) have been already used for the paths~\(P(h_j),j < i.\) Here we use the notation~\(a_n = o(b_n)\) to denote that~\(\frac{a_n}{b_n} \longrightarrow 0\) as~\(n \rightarrow \infty\) and the final relation in~(\ref{e_not_est}) is true since~\(\beta < \frac{1-\alpha}{4}\) strictly, by Theorem statement.

Thus we can always pick a new vertex from~\(S_k(h_i)\) for the relay path~\(P(h_i)\) and continuing iteratively, we get the desired relay RGG~\(G_{\Gamma}.\) The total number of edges in~\(G_{\Gamma}\) is at most~\(\sum_{f}\left(\frac{l(f)}{r_n} + 2\right) = \frac{l_{tot}}{r_n} + 2e_0\) and so combining with~(\ref{e_edge_est}), we get the upper bound in~(\ref{har_str}).

Finally, using~(\ref{eps_bound}) and following a similar analysis as above gives us the upper bound in~(\ref{har_weak}) and this completes the proof of the Theorem.~\(\qed\)


We now equip the edges of~\(G\) with random weights and study relay RGGs with maximum possible weight. Let~\(G_{loc}\) be the random graph as described prior to Theorem~\ref{thm_edge_har} that has~\(v_0+n\) vertices and at most~\(t := v_0 n + {n \choose 2}\) edges. We let~\(Y_{k}, 1 \leq k \leq t\) be positive independent random variables that are also independent of~\(\{X_i\}.\) Using a deterministic rule, we label the edges of the graph~\(G_{loc}\) as~\(f_1,\ldots,f_w, w \leq t\) and for completeness, define~\(f_l, w+1 \leq l \leq t\) to be the straight line segment with endvertices~\(\left(\frac{l}{n^2},0\right)\) and~\(\left(\frac{l+1}{n^2},0\right).\) For each~\(1 \leq k \leq t,\) we then assign weight~\(Y_k\) to the edge~\(f_k.\)

We define the total weight of a relay RGG~\(H = \bigcup_{f \in \Gamma}P(f) \subset G_{loc}\) as~\(W(H) := \sum_{f_k \in H} Y_{k},\)  the sum of weights of  all edges in~\(H.\) We say that~\(H\) is~\(L_n-\)constrained if each~\(P(f)\) has at most~\(L_n\) edges and denote~\(W_n\) to be the maximum weight of an~\(L_n-\)constrained relay RGG, if such a relay RGG exists. Else we set~\(W_n := \sum_{f \in K_n} W(f),\) the total weight of all edges in the complete graph~\(K_n.\) If~\(l_{up}\) is the maximum length of an edge in the deterministic graph~\(\Gamma,\) then~\(L_n \geq \frac{l_{up}}{r_n}\) is a necessary condition for existence of~\(L_n-\)constrained relay RGGs

The following result obtains deviation estimate for~\(W_n\) and also determines sufficient conditions (in terms of~\(L_n\) and~\(r_n\)) for the~\(L^2-\)convergence of~\(W_n,\) appropriately scaled and centred. We recall that~\(l_0\) is the \emph{minimum} length of an edge in~\(\Gamma.\)
\begin{theorem}\label{cost_thm} Suppose~(\ref{gam_cond}) holds with~\(0 < \alpha < 1\) and~\(0 < \beta < \frac{1-\alpha}{2}\) and also suppose that the edge weights are exponentially distributed with unit mean. Also suppose that
\begin{equation}\label{ln_cond}
L_n \geq \frac{l_{up}}{r_n} \text{ and } \frac{r_n^2L_n}{l_0} \longrightarrow 0
\end{equation}
and set~\(\delta_n := e_0 L_n \log{n}.\) For every~\(a > 0,\) there are constants~\(D_1,D_2 >0\) such that
\begin{equation}\label{w_dev_est}
\mathbb{P}\left( D_1\delta_n \leq W_n \leq D_2 \delta_n\right) \geq 1- \frac{1}{n^{1+a}}.
\end{equation}
Moreover, there are constants~\(D_3,D_4>0\) such that~\(D_3 \delta_n \leq \mathbb{E}W_n \leq D_4 \delta_n.\)
\end{theorem}
Thus the maximum weight of an~\(L_n\) constrained relay RGG grows as~\(\delta_n = e_0 L_n\log{n}.\) The term~\(\log{n}\) could be interpreted as the ``gain" obtained in choosing maximum weight paths as opposed to simply selecting a deterministic path of length~\(L_n.\)



\emph{Proof of Theorem~\ref{cost_thm}}: We begin with the upper bound for~\(W_n.\) Let~\(E_{up}\) be the event that~\(Y_{k} \leq M\log{n}\) for each~\(1 \leq k \leq t = v_0n + {n \choose 2},\) where~\(M > 0\) is a constant to be determined later. By~(\ref{gam_cond}), we see that~\(v_0 \leq 2e_0 \leq 2n^{\alpha}\) with~\(\alpha < 1\) and so
\begin{equation}\label{t_est}
t = v_0n + {n\choose 2} \leq 2n^{1+\alpha} + \frac{n^2}{2} \leq n^2.
\end{equation}
Since~\(Y_{k}\) is exponentially distributed with unit mean, we know that\\\(\mathbb{P}(Y_{k} > M\log{n}) = \frac{1}{n^{M}}\) and so by the union bound
\begin{equation}\label{e_up_est}
\mathbb{P}(E_{up}) \geq 1-  \frac{t}{n^{M}}  \geq 1-\frac{1}{n^{M-2}} \longrightarrow 1
\end{equation}
provided~\(M > 2\) strictly. We fix such an~\(M\) and see that if~\(E_{up}\) occurs, then any relay path of length at most~\(L_n\) has weight at most~\(ML_n \log{n}.\) Since there are~\(e_0\) such paths in any relay RGG, we have that~\(W_n \leq e_0ML_n\log{n}\) and this obtains the upper bound for~\(W_n\) in~(\ref{w_dev_est}).

To upper bound~\(\mathbb{E}W_n,\) we use the estimate~\(W_n \leq Y_{tot} := \sum_{k=1}^{t} Y_{k}\) and get that
\begin{equation}\label{ewn_up}
\mathbb{E}W_n \leq e_0ML_n\log{n} + \mathbb{E}Y_{tot}\ind(E^c_{up}),
\end{equation} where~\(\ind(.)\) refers to the indicator function. By the Cauchy-Schwarz inequality, we have that
\begin{equation}\nonumber
\mathbb{E}(Y_{tot} \ind(E^c_{up})) \leq \left(\mathbb{E}Y^2_{tot}\right)^{\frac{1}{2}} \cdot \left(\mathbb{P}(E^c_{up})\right)^{\frac{1}{2}} \leq \left(\mathbb{E}Y^2_{tot}\right)^{\frac{1}{2}} \cdot \frac{1}{n^{M/2-1}}, \nonumber
\end{equation}
using~(\ref{e_up_est}). Moreover the inequality~\(\left(\sum_{i=1}^{k}a_i\right)^2 \leq k \sum_{i=1}^{k}a_i^2\) for positive~\(a_i\) implies that
\begin{equation}\label{y_tot_est_max}
\mathbb{E}Y_{tot}^2 \leq t \sum_{k=1}^{t} \mathbb{E}Y^2_{k} \leq C_1 t^2 \leq 4C_1n^4
\end{equation} for some constant~\(C_1 >0,\)
where the second estimate in~(\ref{y_tot_est_max}) is true by the bounded second moment property of the edge weights (see Theorem statement) and final bound in~(\ref{y_tot_est_max}) is true due to~(\ref{t_est}). Thus
\begin{equation}\label{e_y_tot_est}
\mathbb{E}Y_{tot}\ind(E^c_{up}) \leq \frac{C_2n^{2}}{n^{M/2-1}}
\end{equation} for some constant~\(C_2 >0\) and plugging this into~(\ref{ewn_up}), we get that~\[\mathbb{E}W_n \leq e_0M\log{n} + \frac{C_2n^{2}}{n^{M/2-1}} \leq e_0M\log{n}+1\] provided~\(M\) is sufficiently large. Fixing such an~\(M\) obtains the upper bound for~\(\mathbb{E}W_n.\)



\begin{figure}[tbp]
\centering
\includegraphics[width=5in, trim= 10 340 10 90, clip=true]{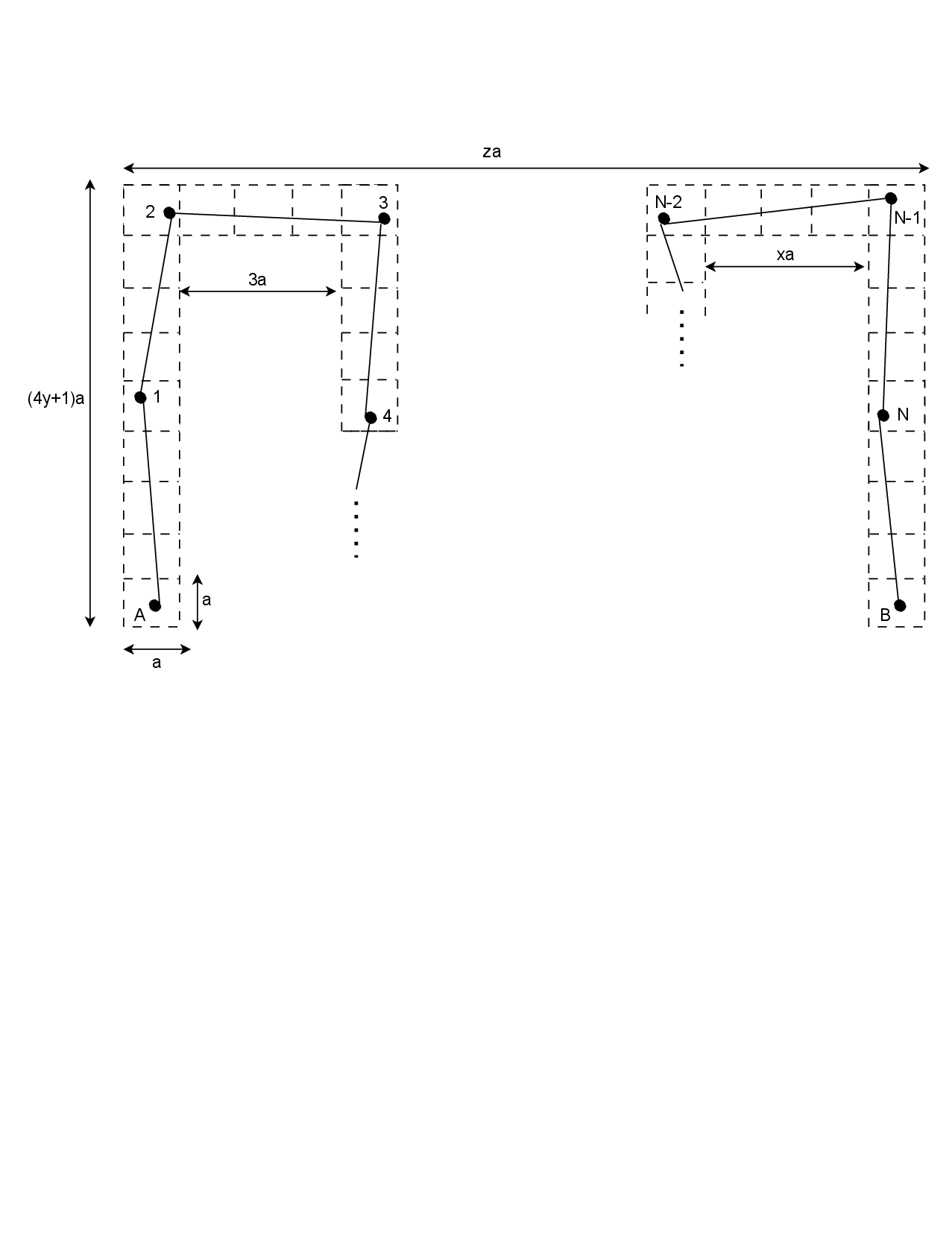}
\caption{The path~\(\tau\) as described in the proof below with~\(U_0 = A,U_{N+1} = B\) and~\(U_i\) denoted as~\(i\) for~\(1 \leq i \leq N.\)}
\label{fig_max_wt}
\end{figure}

Next we prove the lower bound in~(\ref{w_dev_est}) as follows. Initially we assume that~\(\Gamma\) consists of a single edge~\(f = (u,v)\) with~\(u\) as the origin and~\(v = (l(f),0),\) where~\(l(f)\) is the Euclidean length of edge~\(f.\) Let~\(S_0,S_1,\ldots,S_{N},S_{N+1}\) be disjoint~\(a \times a\) squares placed as in Figure~\ref{fig_max_wt} with~\(a = \frac{r_n}{10},\) where~\(S_0\) is the square containing~\(A=u, S_{N+1}\) is the square containing~\(B=v\) and~\(S_i\) is the square containing the dark circle denoted by~\(i,\) for each~\(1 \leq i \leq N.\) Ideally, we would like to choose the squares in such a way that the distance between~\(S_i\) and~\(S_{i+1}\) is always~\(3a\) for any~\(i.\) But because the integer~\(z\) in Figure~\ref{fig_max_wt} need not be of the form~\(4k+1\) with~\(k\) integer,  there is at most one~\(i_0\) for which the distance between~\(S_{i_0}\) and~\(S_{i_0+1}\) between~\(3a\) and~\(7a.\) This is illustrated in Figure~\ref{fig_max_wt}, where the distance between~\(S_{N-2}\) and~\(S_{N-1}\) is~\(xa\) where~\(3 \leq x \leq 7.\)

Based on the discussion in the above paragraph, we get that~\[\frac{z y}{8} \leq \frac{zy}{x+1}\leq   N \leq zy.\] The integer~\(z\) lies between~\(\frac{l(f)}{a}-2\) and~\(\frac{l(f)}{a}\) and since~\(\frac{l_0}{r_n} \rightarrow \infty\) (see~(\ref{gam_cond})), we have that~\(\frac{l(f)}{a}-2 \geq \frac{l(f)}{2a}.\) Setting~\(y = \frac{aL_n}{8l(f)}\) we then get that
\begin{equation}\label{n_est}
\frac{L_n}{16} \leq \frac{l(f)}{2a} \cdot y \leq N \leq \frac{l(f)}{a} \cdot y \leq \frac{L_n}{8} \leq n.
\end{equation}
Moreover since~\(\frac{r_n^2L_n}{l_0} \longrightarrow 0\) by Theorem statement, we see that all the squares~\(\{S_i\}\) are strictly contained in the interior of the unit square~\(S,\) for all~\(n\) large.

Each vertex~\(X_j\) is present in the square~\(S_i\) with probability~\[\int_{S_i} f(y)dy \geq \epsilon_1 x^2  \geq 2C_3 r_n^2\] (see~(\ref{f_eq})) for some constant~\(C_3 > 0\) and so the expected number of vertices of~\(\{X_j\}\) present in~\(S_i\) is at least~\(2C_3nr_n^2.\) If~\(R_i\) is the number of vertices of~\(\{X_j\}\) present in the square~\(S_i,\) then from~(\ref{conc_est_f}),
\begin{equation}\label{ni_est}
\mathbb{P}\left(R_i \geq C_3 nr_n^2\right) \geq 1-e^{-2C_4nr_n^2}
\end{equation}
for some constant~\(C_4 > 0.\) If~\(E_{dense} := \bigcap_{i=1}^{N}\{R_i \geq C_3 nr_n^2\},\) then from~(\ref{ni_est}),~(\ref{n_est}) and the union bound,
\begin{equation}\label{e_dense_est}
\mathbb{P}(E_{dense}) \geq 1- N \cdot e^{-2C_4nr_n^2} \geq 1- n e^{-2C_4nr_n^2} \geq 1-e^{-C_4nr_n^2}
\end{equation}
since~\(nr_n^2\) is much larger than~\(\log{n},\) by Theorem statement.

Suppose~\(E_{dense}\) occurs. We now construct a large weight path using edges of length at least of the order of~\(r_n.\) Indeed, pick the node~\(U_0 = u \in S_0.\) Skip three squares and among all the nodes adjacent to~\(U_0\) in the fourth square~\(S_1,\)  we pick the node~\(U_1\) such that the edge~\((U_0,U_1)\) has the largest possible weight. This is illustrated in Figure~\ref{fig_max_wt}, where~\(S_1\) is the square labelled~\(1.\) We repeat the same procedure with the node~\(U_1\) and the square~\(S_2.\) We continue this procedure until we have chosen the node~\(U_N \in S_N\) and set~\(U_{N+1} = v.\) The length of the path~\(\tau = (U_0,U_1, U_2,\ldots,U_{N},U_{N+1})\) obtained by the selection procedure above is~\(N+2 \geq  \frac{L_n}{16}\) and at most~\(\frac{L_n}{8} < L_n,\) by~(\ref{n_est}). Moreover, by construction, any two vertices in~\(\tau\) are spaced at least~\(\frac{3r_n}{10}\) apart.

To estimate the individual edge weights, we proceed as follows. Since~\(E_{dense}\) occurs, there are at least~\(C_3nr_n^2\) nodes of~\(\{X_j\}\) present in any square~\(S_i.\) Therefore given~\(U_1,\ldots U_j,\) we see that the edge weight~\(V_{j}\) of the edge~\((U_{j},U_{j+1})\) is stochastically dominated from below by~\(\max_{1 \leq j \leq C_3nr_n^2}\zeta_j,\) where~\(\{\zeta_j\}\) are independent exponentially distributed random variables that are also independent of all random variables defined so far. Letting~\({\cal G}_j\) be the sigma-field generated by~\(\{X_i\}_{1 \leq i \leq n} \cup \{U_k\}_{1 \leq k \leq j},\) we have for~\(x > 0\) that
\[\mathbb{P}\left(V_{j} \leq x \mid {\cal G}_j\right) \ind(E_{dense}) \leq (1-e^{-x})^{C_3nr_n^2} \leq (1-e^{-x})^{C_3n^{1-2\beta}}\]
since~\(r_n \geq \frac{1}{n^{\beta}}.\) Setting~\(x = \epsilon \log{n}\) we get that
\begin{equation}
\mathbb{P}\left(V_{j} \leq \epsilon \log{n} \mid {\cal G}_j\right) \ind(E_{dense}) \leq \left(1-\frac{1}{n^{\epsilon}}\right)^{C_3n^{1-2\beta}} \leq \exp\left(-C_3n^{1-2\beta-\epsilon}\right) \nonumber
\end{equation}
and so
\begin{eqnarray}
\mathbb{P}\left(V_{j} \leq \epsilon \log{n}\right) &\leq& \exp\left(-C_3n^{1-2\beta-\epsilon}\right) + \mathbb{P}(E^c_{dense}) \nonumber\\
&\leq& \exp\left(-C_3n^{1-2\beta-\epsilon}\right) + e^{-C_4nr_n^2} \nonumber\\
&\leq& 2\exp\left(-C_3n^{1-2\beta-\epsilon}\right), \nonumber
\end{eqnarray}
by~(\ref{e_dense_est}).

Since~\(\beta < \frac{1}{2},\) we can choose~\(\epsilon > 0\) small so that~\(\beta < \frac{1-\epsilon}{2}\) strictly. Fixing such an~\(\epsilon\) and using the union bound we get
\begin{eqnarray}\label{edge_wt_ax}
\mathbb{P}\left(\bigcap_{j=1}^{N-1} \left\{V_j \geq \epsilon \log{n}\right\}\right) &\geq& 1-2(N-1)\cdot \exp\left(-C_3n^{1-2\beta-\epsilon}\right) \nonumber\\
&\geq& 1-2n\exp\left(-C_3n^{1-2\beta-\epsilon}\right)
\end{eqnarray}
by~(\ref{n_est}). Thus with high probability each edge of~\(\tau\) has weight at least~\(\epsilon \log{n}\) and as argued before, the number of edges in~\(\tau\) is at least~\(\frac{L_n}{16}\) and at most~\(L_n.\) This obtains the lower bound in~(\ref{w_dev_est}) for the case when~\(\Gamma\) has a single edge.

If~\(\Gamma\) has multiple edges, then we order the edges of~\(\Gamma\) as~\(h(1),h(2),\ldots,h(e_0)\) and define an analogous event~\(E_{dense}(i)\) for each~\(1 \leq i \leq e_0.\) Setting~\(E_{dense} := \bigcap_{i=1}^{e_0} E_{dense}(i),\) we get from~(\ref{e_dense_est}) and the union bound that
\[\mathbb{P}(E_{dense}) \geq 1- e_0e^{-C_4nr_n^2}  \geq 1-e^{-C_5nr_n^2}\] for all~\(n\) large and some constant~\(C_5 > 0,\) since~\(e_0\leq n^{\alpha} < n\) by~(\ref{e_not_est}) and~\(nr_n^2\) is much larger than~\(\log{n}\) by Theorem statement.

Assuming that~\(E_{dense}\) occurs, we then iteratively pick~\(P(h_i)\) using the above described procedure in such a way that~\(P(h_i)\) and~\(P(h_j), j < i\) share no vertex in common if~\(h_i\) and~\(h_j\) do not share an endvertex in~\(\Gamma\) and~\(P(h_i)\) and~\(P(h_j)\) share the common endvertex of~\(h_i\) and~\(h_j,\) otherwise. This is possible since any two vertices in~\(P(h_j), j <i\) are spaced at least~\(\frac{3r_n}{10}\) apart in our construction above and the distance between any two vertices within a single square~\(S_i\) is at most~\(\frac{2r_n}{10}.\) Therefore while constructing~\(P(h_i)\) using the~\(a \times a\) squares as described above, we see from~(\ref{ni_est}) that there at least~\(C_3nr_n^2-e_0 \geq C_6nr_n^2\) choices for choosing a vertex in each square, for some constant~\(C_6 >0\) not depending on the choice of~\(i,\) since~\(e_0\) is much smaller than~\(nr_n^2\) by Theorem statement. Following an analogous argument as in the paragraph containing~(\ref{edge_wt_ax}) for the relay RGG construction, we then obtain the lower bound in~(\ref{w_dev_est}). This completes the proof of the Theorem.~\(\qed\)



\begin{thebibliography}{99}
\bibitem{alon} N. Alon and J. Spencer. (2008).
\newblock{\em The Probabilistic Method}.
\newblock{Wiley Interscience}.

\bibitem{david} R. David and U. Feige. (2013).
\newblock{Connectivity of Random High Dimensional Geometric Graphs}.
\newblock{\em International Workshop on Randomization and Approximation Techniques in Computer Science}, pp. 497--512.


\bibitem{diaz} J. D\'iaz, D. Mitsche and X. P\'erez-Gim\'enez. (2002).
\newblock{Large connectivity for dynamic random geometric graphs}.
\newblock{\em IEEE Transactions on Mobile Computing}, \textbf{8}, pp. 821--835.


\bibitem{ferrero} R. Ferrero and F. Gandino. (2017).
\newblock{Analysis of Random Geometric Graph for Wireless Network Configuration}.
\newblock{\em  Tenth International Conference on Mobile Computing and Ubiquitous Network (ICMU), Japan, 2017}, pp. 1-6.

\bibitem{gan_one} G. Ganesan. (2013).
\newblock{Size of Giant Component in Random Geometric Graphs}.
\newblock{\em Annals de l'Institute Henri Poincar\'e}, \textbf{49}, pp. 1130--1140.

\bibitem{gan_two} G. Ganesan. (2018).
\newblock{Stretch and Diameter in Random Geometric Graphs}.
\newblock{\em Algorithmica}, \textbf{80}, pp. 300--330.

\bibitem{goldsmith} A. Goldsmith. (2005).
\newblock{\em Wireless Communications}.
\newblock{Cambridge University Press}.

\bibitem{gupta} P. Gupta and P. R. Kumar. (1998).
\newblock {Critical Power for Asymptotic Connectivity in Wireless Networks}.
\newblock {\em Stochastic Analysis, Control, Optimization and Applications}, pp. 2203--2214.

\bibitem{haenggi} M. Haenggi, J. G. Andrews, F. Baccelli, O. Dousse  and M. Franceschetti. (2009).
\newblock{Stochastic Geometry and Random Graphs for the Analysis and Design of Wireless Networks}.
\newblock{\em  IEEE Journal on Selected Areas in Communications}, \textbf{27}, pp. 1029--1046.

\bibitem{molisch} A. F. Molisch. (2022).
\newblock{\em Wireless Communications: From Fundamentals to Beyond 5G}.
\newblock{Wiley-IEEE Press},~\(3^{rd}\) edition.

\bibitem{penrose} M. Penrose. (2003).
\newblock {\em Random Geometric Graphs}.
\newblock {Oxford University Press}.


\bibitem{pen2} M. Penrose. (2016).
\newblock{Connectivity of Soft Random Geometric Graphs}.
\newblock{\em Annals of Applied Probability}, \textbf{26}, pp. 986--1028.

\bibitem{song} L. Song. (2010).
\newblock{\em Random Graph Models for Wireless Communication Networks}.
\newblock{PhD Thesis, Queen May University, London}

\bibitem{zhang} Z. Zhang, G. Mao and B. D. O. Anderson. (2014).
\newblock{Stochastic Characterization of Information Propagation Process in Vehicular Ad Hoc Networks}.
\newblock{\em IEEE Transactions on Intelligent Transportation Systems}, \textbf{15}, pp. 122--135.


\end{thebibliography}


\subsection*{\em Acknowledgement}
I thank Professors Rahul Roy, Federico Camia and C. R. Subramanian for crucial comments and also thank IMSc and IISER Bhopal for my fellowships.

\subsection*{\em Data Availability Statement}
Data sharing not applicable to this article as no datasets were generated or analysed during the current study.

\subsection*{\em Conflict of Interest}
The authors have no conflicts of interest to declare that are relevant to the content of this article. No funding was received to assist with the preparation of this manuscript.

\end{document}